\documentclass[12pt]{amsart}
\usepackage{amscd,amssymb}

\usepackage[graph,frame,poly,arc]{xy}  
\usepackage[plainpages,backref,urlcolor=blue]{hyperref}

\theoremstyle{plain}
\newtheorem{thm}[subsection]{Theorem}
\newtheorem{lem}[subsection]{Lemma}
\newtheorem{prop}[subsection]{Proposition}
\newtheorem{cor}[subsection]{Corollary}

\theoremstyle{definition}
\newtheorem{rk}[subsection]{Remark}

\newtheorem{ex}[subsection]{Example}

\newtheorem{conjecture}[subsection]{Conjecture}

\numberwithin{equation}{section}
\newcommand{\OO}{{\mathcal O}}

\newcommand{\F}{{\mathcal F}}
\newcommand{\A}{{\mathcal A}}

\newcommand{\HH}{{\mathcal H}}

\newcommand{\B}{{\mathcal B}}
\newcommand{\CC}{{\mathcal C}}

\newcommand{\al}{{\alpha}}
\newcommand{\be}{{\beta}}

\newcommand{\C}{\mathbb{C}}

\newcommand{\PP}{\mathbb{P}}

\DeclareMathOperator{\im}{im}

\DeclareMathOperator{\mult}{mult}
\DeclareMathOperator{\indeg}{indeg}

\begin{document}

\title [On the module of derivations of a line arrangement]
{On the module of derivations of a line arrangement}

\author[Alexandru Dimca]{Alexandru Dimca$^1$}
\address{Universit\'e C\^ ote d'Azur, CNRS, LJAD, France and Simion Stoilow Institute of Mathematics,
P.O. Box 1-764, RO-014700 Bucharest, Romania}
\email{Alexandru.DIMCA@univ-cotedazur.fr}

\thanks{$^1$ partial support from the project ``Singularities and Applications'' - CF 132/31.07.2023 funded by the European Union - NextGenerationEU - through Romania's National Recovery and Resilience Plan.}



\begin{abstract} 
To each multiple point $p$ in a line arrangement $\A$ in the complex projective plane we associate a local derivation $\tilde D_p \in D_0(\A)$.
We show first that these derivations span the graded module of derivations $D_0(\A)$ in all degrees $\geq d -3$, where $d$ is the number of lines in $\A$, see Theorem \ref{thm2} and Theorem \ref{thmG}.  Then, to each local derivation $\tilde D_p \in D_0(\A)$ we associate a polynomial $g_p$ which seems to play a key role in the characterization of the freeness of $\A$, see Theorem \ref{thm4}, as well as in the study of the position of the multiple points of $\A$ with respect to unions of lines, see Corollary \ref{cor101} and Conjecture \ref{conj10}. Corollary \ref{cor1G} gives a result of an independent interest, namely a lower bound for the maximal exponent of a plane curve having a line as an irreducible component.

\end{abstract}

\maketitle

%\tableofcontents

\section{Introduction} 

Let $S=\C[x,y,z]$ be the graded polynomial ring in the variables
$x,y,z$ with complex coefficients and let $f \in S_d$ be a homogeneous polynomial of degree $d$. Let $f_x, f_y$ and $f_z$ denote the partial derivatives of $f$ with respect to $x, y$ and $z$ respectively. Recall that a Jacobian syzygy  for $f$ is a triple
 \begin{equation}
\label{eq1}
\rho=(a,b,c) \in S^{3} \text{ such that } af_x+bf_y+cf_z=0.
\end{equation}
With such a syzygy $\rho$, one can associate the following derivation 
 \begin{equation}
\label{eq11}
D_{\rho}=a{\partial_x}+ b\partial_y+ c\partial_z \in Der(S),
\end{equation}
of the polynomial ring $S$, such that $D_{\rho}(f)=0$, i.e. the derivation $D_{\rho}$ kills the polynomial $f$.
Consider the graded $S$-module of Jacobian syzygies of $f$, namely
$$Syz(f)=\{ \rho=(a,b,c) \in S^{3} \text{ such that } af_x+bf_y+cf_z=0.   \}.$$
Note that $Syz(f)=D_0(f)$, the $S$-module of all derivations killing $f$, as defined in \eqref{eq11}. {\it  In the sequel we often identify the syzygy $\rho$ with the corresponding derivation $D_{\rho}$.}
We say that the curve $C:f=0$ is an {\it $s$-syzygy curve} if  the module $D_0(f)$ is minimally generated by $m$ homogeneous syzygies, say $\rho_1,\rho_2,...,\rho_s$, of degrees $d_j=\deg \rho_j$ ordered such that $$d_1\leq d_2 \leq \ldots \leq d_s.$$ 
We call these degrees $(d_1, \ldots, d_s)$ the {\it exponents} of the curve $C$. 
The smallest degree $d_1$ is sometimes denoted by ${\rm mdr}(f)$ and is called the minimal degree of a Jacobian relation for $f$. 

In this paper we consider the case when $C$ is a line arrangement $\A:f=0$, and try to detect properties of the module $D_0(f)$ of derivations killing $f$ which depend only on the combinatorics of the arrangement $\A$.  Our first main result is the following.

\begin{thm}
\label{thmD} For the arrangements $\A:f=0$ of $d$ lines in $\PP^2$ one has the following.

\begin{enumerate} 
\item The dimension $\dim D_0(f)_j$ is combinatorially determined for any $j \geq d-3$. More precisely, for $j \geq d-3$, one has
$$\dim D_0(f)_{j}=
\tau(\A)+3 \binom{j+2}{2}-\binom{d+j+1}{2}.$$
In particular
$$\dim D_0(f)_{d-3}= \tau(\A)-\binom{d}{2} \text{ and } \dim D_0(f)_{d-2}= \tau(\A)-\binom{d-1}{2}.$$

\item The dimension $\dim D_0(f)_{d-4}$ is not combinatorially determined in general.

\end{enumerate} 

\end{thm}
As a direct consequence of Theorem \ref{thmD} we have the following.

\begin{cor}
\label{corD} For any arrangement $\A:f=0$ of $d$ lines in $\PP^2$ one has the following, where $m_p$ is the multiplicity of the point $p \in \A$ and the sum is over all points of multiplicity at least $3$.

\begin{enumerate} 
\item
$$\dim D_0(f)_{d-3}= \sum_p \binom{m_p-1}{2}.$$
 In particular $\dim D_0(f)_{d-3}=0$ if and only if $\A$ has only double points.

\item $$\dim D_0(f)_{d-2}= \sum_p \binom{m_p-1}{2} +d-1.$$
 In particular $\dim D_0(f)_{d-2}=d-1$ if and only if $\A$ has only double points.
\end{enumerate}  
 
\end{cor}
A second theme of this paper is to construct some derivations in $D_0(\A)$ starting from the combinatorics of $\A$. 
More precisely, let $\A:f=0$ be a line arrangement of $d$ lines in $\PP^2$ and let 
$$p=(u:v:w)$$ be a point in $\A$ of multiplicity $m_p \geq 2$. We assume in this paper that $\A$ is not a pencil of lines, that is we assume that $m_p<d$.
Then the defining equation $f=0$ can be written as
$$f=f_{1p}f_{2p}=0,$$
where $f_{1p}=0$ (resp. $f_{2p}=0$) is the equation of the subarrangement $\A_1$ (resp. $\A_2$) of $\A$ formed by all the lines in $\A$ passing (resp. not passing) through the point $p$. The following result tells that to the multiple point $p$ one can associate in a rather natural way a derivation in $D_0(\A)_{d-m_p}$. If $p$ is any point in the complex projective plane $\PP^2$, 
 and $D_{\rho} $ is a derivation as in \eqref{eq11}, then we set
 \begin{equation}
\label{eq2}
D_{\rho}(p)=a(p)\partial_x+ b(p)\partial_y+ c(p)\partial_z ,
\end{equation}
that is we evaluate all the coefficients of $D_{\rho}$ at $p$, and get a constant derivation, well defined up-to a nonzero constant factor.
\begin{thm}
\label{thm1} For any arrangement $\A:f=0$ of $d$ lines in $\PP^2$, and
with the above notation, we consider the following two derivations
$$D_p=u\partial_x+ v\partial_y+ w\partial_z$$
and
$$\tilde D_p=f_{2p}D_p- \frac{D_p(f_{2p})}{d}E,$$
where $E=x{\partial _x}+ y{\partial_ y}+ z{\partial _z}$ is the Euler derivation.
Then the following hold.

\begin{enumerate} 
\item $D_p \in D_0(f_{1p})$ and $\deg D_p=0$.
\item $\tilde D_p \in D_0(f)$ and $\deg \tilde D_p=d-m_p$.

\item $\tilde D_p(p) \ne 0$.

\item Let $q\ne p$ be another multiple point of $\A$. Then the coefficients of $\tilde D_p$ vanish at $q$ of order at least $m_q-2$.

\item For any point $q \in \PP^2$, one has  $\tilde D_p(q) =0$ if and only if $q$ is a multiple point of the arrangement $\A_2:f_{2p}=0$.

\end{enumerate} 

\end{thm}
Such local derivations $\tilde D_p$ were already considered in our papers
\cite[proof of Theorem 1.2]{Mich} and \cite[proof of Lemma 3.2]{minT}, but now we consider them in more details. The vector spaces spanned by the multiples of such local derivations is quite large in $D_0(f)$. For instance, we have the following result.

\begin{thm}
\label{thm2}
For any arrangement $\A:f=0$ of $d$ lines  
the vector space $D_0(f)_{d-3}$ is the direct sum
$$\oplus_pS_{m_p-3}\tilde D_{p},$$
where $p$ runs through the set of points $p \in \A$ of multiplicity $m_p \geq 3$.
\end{thm}
The following result is a direct consequence of Theorem \ref{thm2}, to be compared with Theorem \ref{thmD} $(2)$.
\begin{cor}
\label{cor20}
For any arrangement $\A:f=0$ of $d$ lines in $\PP^2$ and for any $j \geq 4$, one has the following inequality
$$\dim D_0(f)_{d-j} \geq \sum _p \binom{m_p-j+2}{2},$$
where the sum is over all multiple points $p \in \A$ such that $m_p \geq j$.
\end{cor}
The example of arrangement in the proof of Theorem \ref{thmD} $(2)$
shows that the inequality in Corollary \ref{cor20} may be strict.
As a complement to Theorem \ref{thm2}, we can prove the  following key result.

\begin{thm}
\label{thmG}
For any arrangement $\A$ of $d$ lines in $\PP^2$ which is not a pencil, and 
for any $j \geq d-2$, the local derivations $\tilde D_p$ generate
$D_0(\A)_j$, when $p$ runs through all the multiple points of $\A$ of multiplicity $ m_p \geq 2$. More precisely, one has
$$D_0(\A)_j=\sum_p S_{j+m_p-d}\tilde D_p,$$
for any $j \geq d-2$, where this sum is not a direct sum in general.
\end{thm}
\begin{rk}
\label{rkG}
Since all the generators of the graded $S$-module $D_0(\A)$ have degrees $\leq d-2$, see for instance \cite{Sch}, it follows that to prove Theorem \ref{thmG} it is enough to prove the case $j=d-2$.
\end{rk}
The main ingredient in the proof of Theorem \ref{thmG} is the following result, which we state and prove in a more general setting than needed in this paper.
Let $C=C' \cup L$ be a reduced curve, where $C':f'=0$ is the curve obtained by deleting the line $L$ from $C$. 
Let $\A'':f''=0$ be the multiarrangement on $L$ defined by the restriction of $f'$ to $L$.
We choose the coordinates such that $L:x=0$ and, then $f=xf'$ and $f''=\overline f'$,
where we define for any polynomial $g \in S$
$$\overline g(y,z)=g(0,y,z).$$
In this setting there is an exact sequence of graded $S$-modules
 \begin{equation}
\label{eq0G}
0 \to D(f')(-1) \stackrel{u} \longrightarrow D(f) \stackrel{v} \longrightarrow D(f''),
\end{equation}
where the morphisms $u$ and $v$ are given by
\begin{equation}
\label{eq01G}
u(\delta)=x\delta \text{ and } v(a,b,c)=(\overline b,\overline c)
\end{equation}
and $D(f'')$ is regarded as an $S$-module via the morphism $S \to \C[y,z]$, $g \mapsto \overline g$,
see \cite[Proposition 4.45]{OT} in the case of line arrangements  and \cite{STY} in the general case.
One may ask whether a similar exact sequence exists for the modules
$D_0(f')$, $D_0(f)$ and $D_0(f'')$, namely
\begin{equation}
\label{eq2G2}
0 \to D_0(f')(-1) \stackrel{u_0} \longrightarrow D_0(f) \stackrel{v_0} \longrightarrow D_0(f'').
\end{equation}
This answer is known to be positive and can be recovered from  \cite{Sch} or \cite{STY}, where $D_0(f')$, $D_0(f)$ and $D_0(f'')$ are identified to $D(f')/S\cdot E$, $D(f)/S\cdot E$ and respectively to $D(f'')/S\cdot \overline E$. It seems that this identification gives simpler, but less precise formulas for $u_0$ and $v_0$ compared to our next result.
\begin{thm}
\label{thm1G} 
With the above notation, the sequence in \eqref{eq2G2} is exact if we define the morphisms $u_0$ and $v_0$ of graded $S$-modules as follows.
\begin{enumerate} 

\item For $\delta=(a,b,c)\in D_0(f')$, we set
$$u_0(\delta)=x \delta-\frac{aE}{d}$$
with $d= \deg f$ and $E=(x,y,z)$ the Euler derivation, 

\item Any $\delta \in D_0(f)$ has the form $\delta=(xa_1,b,c)$ with $a_1 \in S$, and with this notation we set
$$v_0(\delta)=(\overline b,\overline c)+\frac{\overline a_1}{d-1}(y,z).$$

\end{enumerate} 

\end{thm}
A nice consequence of Theorem \ref{thm1G} and of \cite{POG} is  the following.

\begin{cor}
\label{cor1G} 
Consider a reduced plane curve $C=C' \cup L$ in $\PP^2$, where $C':f'=0$ is another reduced curve and $L$ is a line not contained in $C'$.
Let $r$ be the number of points in the set given by the intersection  $C' \cap L$. Then the following hold.

\begin{enumerate} 

\item The morphism $u_0$ induces isomorphisms 
$$D_0(f')_{k-1} \to D_0(f)_k$$
 for all $k\leq r-2$.

\item For $k$ large enough one has 
$$\dim D_0(f')_{k-1} < \dim D_0(f)_k.$$
 In particular,
the maximal exponent $d_s$ of $C$ satisfies
$$d_s \geq r-1.$$

\end{enumerate} 

 \end{cor}

Let $\A:f=0$ be a line arrangement and fix a minimal degree non-zero syzygy $\rho \in D_0(f)$, say $\rho=(a,b,c)$. For any homogeneous syzygy $\rho'=(a',b',c') \in D_0(f)$, consider the determinant $\Delta(\rho')=\det M(\rho')$ of the $3 \times 3$ matrix $M(\rho')$ which has as first row $x,y,z$, as second row $a,b,c$ and as third row $a',b',c'$. Then it turns out that $\Delta(\rho')$ is divisible by $f$, see \cite{Dmax}. For any multiple point $p \in \A$ we consider the corresponding local derivation $\tilde D_p$ and define a homogeneous polynomial $g_p \in S$ with $\deg g_p=d_1+1-m_p$ by setting
\begin{equation} \label{B44}
 g_p=\frac {\Delta(\tilde D_p)}{f}.
\end{equation}
These polynomials $g_p$ enter in the following description of free line arrangements.

\begin{thm}
\label{thm4}
With the above notation, let $m$ be the maximal multiplicity of a point in $\A$. Then the following hold.
\begin{enumerate} 
\item If $m > d/2$, then $d_1=d-m$ and the line arrangement $\A$ is free if and only if
$$\tau(\A)=(d-1)^2-(d-m)(m-1).$$

\item If $m \leq d/2$, then either $d_1=m-1$ and the arrangement $\A$ is free with exponents $d_1 <d_2=d-m$, or
$$m \leq d_1  \leq d-m.$$

\item If $d-2 \geq m \geq 2$, then
the line arrangement $\A$ is free if and only if the ideal
$I(\rho) \subset S$ spanned by all the polynomials $g_p$ with $p\in \A$ satisfying 
\begin{equation} \label{E1}
3 \leq m_p \leq \frac{d+1}{2}
\end{equation}
 has an empty zero set in $\PP^2$.

\end{enumerate} 
\end{thm}

Note that {\it the polynomial $g_p$ depends to some extent on the combinatorics}, but not entirely, since we cannot describe $\rho$ in terms of the combinatorics,
even when $\A$ is a free line arrangement. However, there is one situation where this can be done, and this leads to the following  construction.
Let $\A$ be a line arrangement in $\PP^2$, and fix a point $q \in \A$, such that $m_q=m$ is the maximal multiplicity of points in $\A$. For any other point $p \ne q$, $p \in \A$ with $m_p \geq 3$, we define a line arrangement
$$\B_p:h_p=0$$
 as follows. If the points $p$ and $q$ are connected by a line $L \in \A$, then $\B_p$ consists simply of all the lines in $\A$ which contain neither $p$ nor $q$.
If the points $p$ and $q$ are not connected by a line $L \in \A$, then $\B_p$ consists of all lines in $\A$ which contain neither $p$ nor $q$, plus the line $L(p,q)$ determined by the points $p$ and $q$.
Note that in both cases one has 
$$|\B_p|=d-m_p-m_q+1.$$

\begin{prop}
\label{prop40}
With the above notation, if the condition
$$m>\frac{d}{2}$$
is satisfied, then one may choose $\tilde D_q$ as the minimal degree derivation  $\rho \in D_0(f)$ and with this choice one has
$$g_p=h_p$$
up to a nonzero constant factor, for any point $p$ distinct from $q$.
\end{prop}

The polynomials $g_p$ introduced above enter also in the following result which gives in part (2) a detailed geometrical description of line arrangements with $m(\A)=d_1$.

\begin{thm}
\label{thm1000}
Let $\A:f=0$ be a line arrangement having the first exponent $d_1 \geq 1$ and consider a point $p \in \A$ of multiplicity $m_p$ such that
$$m_p \leq d_1+1 <d-m_p+1.$$
Then all the other multiple points of $\A$ which are not connected to $p$ by lines in $\A$ are situated on a curve $\CC_p$ of degree at most $d_1+1-m_p$, which is a union of some irreducible components of the curve $g_p=0$.
In particular, one has the following.

\begin{enumerate} 

\item If $m_p=d_1+1$, then $p$ is a modular point and the line arrangement $\A$ is supersolvable.

\item If $m_p=d_1$ and $m_p=m(\A)$, the maximal multiplicity of points in $\A$, then all the other multiple points $q \in \A$, $q \ne p$, which are not connected to $p$ by lines in $\A$ are situated on a line $L$.
If $Q$ denotes the set of such points $q$, then one has
$$ d-2m_p \leq \sum_{q\in Q} (m_q-1)$$
and the equality holds if and only if $\A$ is free with exponents $(m_p,d-m_p-1)$.
Moreover, if this inequality is strict, then $\A$ is a plus-one generated arrangement with exponents $(m_p,d-m_p, d_3)$, where
$$d_3=m_p-1+ \sum_{q\in Q} (m_q-1).$$

\end{enumerate} 

\end{thm}
An alternative proof of the claim $(1)$ in Theorem \ref{thm1000} can be obtained using \cite[Theorem 1.2]{Mich} plus \cite[Proposition 2.3]{UNEXP}. For the claim $(2)$ we note that there are examples when the line $L$ cannot be chosen to be in $\A$, see Example \ref{ex10}.
\begin{cor}
\label{cor101}
Let $\A:f=0$ be a line arrangement having the first exponent $d_1 \geq 1$. 
Then all the multiple points of $\A$  are situated on a curve $\CC$ of degree at most $d_1+1$.
\end{cor}
In fact, this result holds in a more general situation, when the line arrangement is replaced by an arrangement of smooth curves meeting transversally, see \cite[Proposition 3.6]{To} and \cite[Remark 1.3]{BT}.
Our proof shows that the curve $\CC$ contains at least $m_p$ lines, so one may hope that the following more precise statement holds.
\begin{conjecture}
\label{conj10}
Let $\A:f=0$ be a line arrangement in $\PP^2$ having the first exponent $d_1 \geq 1$. 
Let $N$ be the minimal number of lines  in $\PP^2$, such that all the multiple points of $\A$  are situated on a union of $N$ lines.
Then 
$$N \leq d_1+1.$$
\end{conjecture}
As explained in the proof of Corollary \ref{cor101}, Conjecture \ref{conj10} is equivalent to asking that the polynomial $g_p$ is a product of linear factors for at least one multiple point $p$ of any arrangement $\A$.
Note that this happens for any $p \ne q$ in the conditions of Proposition \ref{prop40} above, case which is discussed again in Theorem \ref{thm100} $(3)$ below without the use of polynomials $g_p$.
The next two results support Conjecture \ref{conj10}. The first one is a direct consequence of Theorem \ref{thm1000}, $(2)$.
\begin{cor}
\label{cor1000}
Let $\A:f=0$ be a line arrangement having the first exponent $d_1 \geq 1$ and  a point $p \in \A$ with multiplicity $m_p =d_1$. Then Conjecture \ref{conj10} holds for the arrangement $\A$.
\end{cor}
We have also the following.
\begin{thm}
\label{thm100}
Suppose that the line arrangement $\A$ has a point $p$ of multiplicity $m \geq 2$ such that $d_1=m-1$ or $d_1=d-m$, where $d= |\A|$.
Then Conjecture \ref{conj10} holds for the arrangement $\A$ and moreover the $N$ lines containing all the multiple points of $\A$ can be chosen in $\A$. In particular, this property holds
for the following types of arrangements.
\begin{enumerate} 

\item $\A$ has only double points and $d \geq 3$.

\item $\A$ has only double points except for $n_3$ triple points with $1 \leq n_3 \leq 3$ and $d \geq 4$.

\item $\A$ has a point $p$ of multiplicity $m$ such that $2m \geq |\A|$.

\item $\A$ is  supersolvable.

\end{enumerate} 

\end{thm}

Here is the organization of this paper. In Section 2 we recall some basic facts and then we prove Theorem \ref{thmD} and  Corollary \ref{corD}.
In Section 3 we give the proofs of Theorems \ref{thm1} and \ref{thm2}. We also state and prove a related result involving $D_0(f)_{d-2}$, see Theorem \ref{thm3} and Example \ref{ex2}.

In Section 4 we explain how the local derivations $\tilde D_p$  behave under addition-deletion mappings, and prove Theorem \ref{thmG}.
In Section 5 we show first that the polynomial $g_p$ from \eqref{B44} is nonzero when $\tilde D_p$ is not a constant multiple of minimal degree
syzygy $\rho$, see Proposition \ref{prop4}. Then we prove Theorem \ref{thm4} by using a key property of Bourbaki ideals associated with plane curves. Proposition \ref{prop40} is also proved here. 
In Section 6  we prove Theorems \ref{thm1000} and \ref{thm100} as well as Corollary \ref{cor101}

\medskip

We would like to thank Takuro Abe,  Dan Bath and \c Stefan Toh\u aneanu for some useful discussions on the results in this paper.

\section{The proof of Theorem \ref{thmD} and of Corollary \ref{corD}} 

We denote by $J_f$ the Jacobian ideal of $f$, i.e. the homogeneous ideal in $S$ spanned by $f_x,f_y,f_z$, and  by $M(f)=S/J_f$ the corresponding graded ring, called the Jacobian (or Milnor) algebra of $f$.
 Let $I_f$ denote the saturation of the ideal $J_f$ with respect to the maximal ideal ${\bf m}=(x,y,z)$ in $S$ and consider the local cohomology group 
 $$N(f)=I_f/J_f=H^0_{\bf m}(M(f)).$$
 The graded $S$-module  $N(f)$ satisfies a Lefschetz type property with respect to multiplication by generic linear forms, see  \cite{DPop}. This implies in particular the inequalities
\begin{equation}
\label{eqDP1} 
0 \leq n(f)_0 \leq ...\leq n(f)_{\lceil  T/2 \rceil} \geq n(f)_{\lceil T/2\rceil+1} \geq ...\geq n(f)_T \geq 0,
\end{equation}
where $T=3d-6$, $n(f)_k=\dim N(f)_k$ for any integer $k$, and for any real number $u$,  $\lceil     u    \rceil $ denotes the round up of $u$, namely the smallest integer $U$ such that $U \geq u$.
The self duality of the graded $S$-module $N(f)$, see \cite{Se}, implies that
\begin{equation}
\label{eqDP2} 
 n(f)_{k} = n(f)_{T-k}, 
\end{equation}
for any integer $k$.
The general form of the minimal resolution for the Milnor algebra $M(f)$ of such a curve $C:f=0$ 
is 
\begin{equation}\label{RMM}0\to \displaystyle{\oplus_{i=1}^{s-2}S(-e_i)}\to \displaystyle{\oplus_{i=1}^{s}S(1-d-d_i)}\to S^3(1-d)\to S,\end{equation} 
with $e_1\leq e_2\leq\ldots\leq e_{s-2}$ and $1 \leq d_1\leq d_2\leq\cdots\leq d_s$. It follows from \cite[Lemma 1.1]{HS12} that one has $$e_j=d+d_{j+2}-1+\epsilon_j,$$ for $j=1,\ldots, s-2$ and some integers $\epsilon_j\geq 1.$
The minimal resolution of $N(f)$ obtained from (\ref{RMM}), by \cite[Proposition 1.3]{HS12}, is
 $$0\to \displaystyle{\oplus_{i=1}^{s-2}S(-e_i)}\to\displaystyle{\oplus_{i=1}^sS(-\ell_i)}\to\displaystyle{\oplus_{i=1}^s}S(d_i-2(d-1))\to \displaystyle{\oplus_{i=1}^{s-2}S(e_i-3(d-1))},$$
 where $\ell_i=d+d_i-1$. It follows that 
\begin{equation}\label{sigma} 
 \sigma(C)=3(d-1)-e_{s-2}=2(d-1)-d_s-\epsilon_{s-2},
\end{equation} 
where $$\sigma(C)=\min\left\{j : n(f)_j\neq 0\right\}=\indeg(N(f)).$$

\subsection{The proof of Theorem \ref{thmD}} 

Using \cite[Equation (2.3)]{Dmax}, we get that for any integer $k$ one has
$$\dim D_0(f)_{k+1}-\dim N(f)_{d+k}+\dim D_0(f)_{d-5-k}=$$
$$=\tau(\A)+3 \binom{k+3}{2}-\binom{d+k+2}{2}.$$
If $k\geq d-4$, then $D_0(f)_{d-5-k}=0$ since $d-5-k<0$.
Next we show that 
\begin{equation}
\label{eqNN} 
 N(f)_{2d-4}=0,
\end{equation}
or, equivalently, $N(f)_{d-2}=0$, recall \eqref{eqDP2}. In view of
 \eqref{eqDP2}, this is equivalent to showing
 $$\sigma(\A) \geq d-1.$$
 Note that the minimal resolution \eqref{RMM} gives rise to the following
 minimal resolution for $D_0(f)$:
 \begin{equation}\label{RMA}0\to \displaystyle{\oplus_{j=1}^{s-2}S(-(d_{j+2}+\epsilon _j))}\to \displaystyle{\oplus_{i=1}^{s}S(-d_i)}.\end{equation} 
 It is known that the Castelnuovo-Mumford regularity of the $S$-module
 $D_0(\A)$ is at most $d-2$,  see \cite[Corollary 3.5]{Sch} and also
 \cite{CM}. This implies that
 \begin{equation}
\label{eqNN3} 
 d_s \leq d-2 \text{ and } d_s+ \epsilon_{s-2} -1 \leq d-2.
\end{equation} 
Using the second inequality we get
$$ \sigma(\A)=2(d-1)-d_s-\epsilon_{s-2} \geq d-1,$$
which is exactly what we need to get  \eqref{eqNN}.
Now, for $k \geq d-4$, we have 
 $d+k \geq 2d-4$ and hence $N(f)_{d+k}=0$. Indeed,  we have shown that $N(f)_{2d-4}=0$ and the claim follows using \eqref{eqDP1} since
$$T/2=3(d-2)/2  \leq 2d-4$$
for $d\geq 2$. This completes the proof of claim $(1)$. 

To prove $(2)$, recall the Ziegler pair of two arrangements $\A:f=0$ and $\A':f'=0$ of $d=9$ lines, such that $D_0(f)_5 \ne 0$ and $D_0(f')_5=0$,
see \cite{Zi, ZP}.
\endproof

\subsection{The proof of Corollary \ref{corD}} 
Use Theorem \ref{thmD} and the following two well known formulas
$$\tau(\A)=\sum_p (m_p-1)^2 \text{ and } \sum_p \binom{m_p}{2}=\binom{d}{2}.$$
\endproof

\section{The proofs of Theorems \ref{thm1}  and \ref{thm2} and a related result}

\subsection{The proof of Theorem \ref{thm1}}

To prove $(1)$, it is enough that for any line $L:\ell=0$ in $\PP^2$, one has that $p \in L$ if and only if $D_p(\ell)=0$. Indeed, $f_{1p}$ is just a product of such linear forms $\ell$ with $D_p(\ell)=0$.

One can prove $(2)$ either by a direct computation or using
\cite[Theorem 5.1]{DIS}.

To prove $(3)$, we compute the coefficients $A,B,C$ of ${\partial _x}$, ${\partial_y}$ and respectively ${\partial_z}$ in $\tilde D_p(p)$ and we get
$$(A:B:C)=\frac{m_p}{d}f_{2p}(p)(u:v:w),$$
since $ D_p(f_{2p})(p)=d_2f_{2p}(p)$, with $d_2= \deg f_{2p}=d-m_p >0$.

To prove $(4)$, note that $f_{2p}$ vanishes at $q$ of order $m_q$ (resp. $m_q-1$) if the line $L(p,q)$ determined by $p$ and $q$ is not in $\A$ (resp. if $L(p,q)$ is in $\A$). Hence $f_{2p}$ and $D_p(f_{2p})$ vanish at $q$ of order at least $m_q-2$.

Finally, to prove  $(5)$, note that we may suppose $p=(1:0:0)$. With 
this choice one has
$$\tilde D_p=\left(f_{2p}-\frac{xf_{2p,x}}{d}\right)\partial_x-\frac{yf_{2p,x}}{d}\partial_y-\frac{zf_{2p,x}}{d}\partial_z.$$
A point $q=(u:v:w) \ne p$ has either $v \ne 0$ or $w \ne 0$. It follows that $\tilde D_p(q)=0$ if and only if $f_{2p,x}(q)=f_{2p}(q)=0$, in particular  $q \in \A_2$. If
 $q$ is a point of multiplicity at least two in $\A_2$, then clearly $\tilde D_p(q) =0$, since $f_{2p}$ and all of its first derivatives 
 vanish at $q$.
 If $q$ is a simple point of $\A_2$, and $L:\ell=0$ is the unique line in $\A_2$ containing $q$, then  one has
$$( f_{2p,x}(q):f_{2p,y}(q): f_{2p,z}(q))=(\al:\be:\gamma),$$
where $\ell=\al x+ \be y + \gamma z$. Indeed, both sides give the coefficients of the tangent line $L:\ell=0$ to $\A_2$ at the simple point $q$. It follows that $\al =0$ and hence $L$ is a line containing $p$, a contradiction, since $L$ was supposed to be a line in $\A_2$.
Hence this second case is impossible, and this proves the claim $(5)$.

\subsection{The proof of Theorem \ref{thm2}} 

We show that there is no notrivial relation
$$\sum_pP_p\tilde D_p=0,$$
where $P_p \in S_{m_p-3}$ are not all zero.
Let $M$ be a minimal subset of the set of points in $\A$ of multiplicity at least $3$ such that 
$$\sum_{p\in M}P_p\tilde D_p=0,$$
where $P_p \in S_{m_p-3}$ are  nozero for any $p \in M$.
Let $q \in M$ be a point of minimal multiplicity, that is $m_q \leq m_p$ for all $p \in M$.
Then using Theorem \ref{thm1}, claims $(3)$ and $(4)$, we see that the coefficient $P_q$ has to vanish at $q$ of order at least $m_q-2$.
Indeed, all the other coefficients of $\tilde D_p$ for $p \ne q$ have this property, and at least one coefficient of $\tilde D_q$ does not vanish at $q$. Since $P_q$ is a homogeneous polynomial of degree $m_q-3$,
this implies $P_q=0$ a contradiction.

We have also the following result.

\begin{thm}
\label{thm3}
Let $p,p'$ and $p''$ be three non-collinear multiple points in an arrangement $\A$ of $d$ lines, such that at most two of the three lines determined by $p$, $p'$ and $p''$ are in $\A$. Then 
$D_0(f)_{d-2}$ contains the direct sum
$$S_{m-2}\tilde D_p \oplus S_{m'-2}\tilde D_{p'} \oplus S_{m''-2}\tilde D_{p''} ,$$
where $m=\mult_p(\A) \geq 3 $, $m'=\mult_{p'}(\A) \geq 3 $ and $m''=\mult_{p''}(\A) \geq 3 $.
\end{thm}

\proof

Since $p,p'$ and $p''$  are non-collinear, we may choose the coordinates such that
$$p=(1:0:0), \  p'=(0:1:0) \text{ and } p''=(0:0:1).$$
With this choice, we get 3 derivatives as in Theorem \ref{thm1}, namely
$$\tilde D_p=f_{2p}\partial_x- \frac{f_{2p,x}}{d}E, \ \tilde D_{p'}=f_{2p'}\partial_y- \frac{f_{2p',y}}{d}E \text{ and } \tilde D_{p''}=f_{2p''}\partial_z- \frac{f_{2p'',z}}{d}E.  $$
Let $P,P'$ and $P''$ be homogeneous polynomials in $S$ of degree
$m-2, \ m'-2$ and respectively $m''-2$. We show that 
$$\Sigma =P\tilde D_p+P'\tilde D_{p'}+P''\tilde D_{p''}=0$$
implies $P=P'=P''=0$. By symmetry, we can assume that the line $x=0$ is not in $\A$.
Let $A$ be the coefficient of $\partial _x$ in the sum $\Sigma \in D_0(f)_{d-2}$. We get
$$dA=P(df_{2p}-xf_{2p,x})-P'xf_{2p',y}-P''xf_{2p'',z}.$$
Since $x=0$ is not in $\A$, it follows that $f_{2p}$ is not divisible by $x$.
Then $A=0$ yields $P=xP_1$, for some homogeneous polynomial $P_1$ of degree $m-3$.
Note that $\deg f_{2p}=d-m$ and $f_{2p}(p) \ne 0$ implies that
$x^{d-m}$ occurs in $f_{2p}$ with a non-zero coefficient. It follows that
the factor $(df_{2p}-xf_{2p,x})$ does not vanish at $p$. Now note that
$f_{2p'}$ and $f_{2p''}$ both vanish at $p$ of order at least $m-1$.
More precisely, $f_{2p'}$ vanishes at $p$ of order  $m-1$ if the line
$z=0$ is in $\A$, and of order $m$ if the line
$z=0$ is not in $\A$.
Hence
$A=0$ would imply that $P_1$ vanishes at $p$ of order at least $m-2$.
Since $\deg P=m-3$, this is possible only if $P_1=0$. Hence $P=0$, and we are left with a relation
$$\Sigma_1 =P'\tilde D_{p'}+P''\tilde D_{p''}=0.$$
Moreover $A=0$ and $P=0$ give
$$P'f_{2p',y}+P''f_{2p'',z}=0$$
Let $B$  be the coefficient of $\partial _y$  in the sum $\Sigma_1$.
Then we get
$$dB=P'(df_{2p'}-yf_{2p',y})-P''yf_{2p'',z}=0=dP'f_{2p'}-y(P'f_{2p',y}+P''f_{2p'',z})=dP'f_{2p'}=0.$$
This implies $P'=0$ and then of course $P''=0$ as well. This ends the proof of our claim.

\endproof

\begin{ex}
\label{ex2}
The vector space $D_0(f)_{d-2}$ does not have a the direct sum decomposition similar to that in Theorem \ref{thm2} and the extra assumption in Theorem \ref{thm3} is necessary.
Consider for instance the line arrangement
$$\A: xyz(x-y)(y-z)(x-z)=0.$$
We consider the triple points of $\A$, that is the points
$$q=(1:1:1), \ p=(1:0:0), \  p'=(0:1:0) \text{ and } p''=(0:0:1).$$
It is known that $\A$ is free with exponents $(2,3)$ and hence
$$\dim D_0(f)_4= \dim S_2+\dim S_1=6+3=9.$$
On the other hand 
$$ S_1 \tilde D_p+S_1 \tilde D_{p'}+S_1  \tilde D_{p''}  \ne D_0(f)_4,$$ 
since $\tilde D_p(q)=\tilde D_{p'}(q)=\tilde D_{p''}(q)=0$ by Theorem \ref{thm1} $(4)$ and $D=x\tilde D_q \in D_0(f)_4$ satisfies $D(q) \ne 0$ by Theorem 
\ref{thm1} $(3)$. It follows that there are non-trivial relations
$\Sigma=0$ as in the proof of Theorem \ref{thm2}, with $\deg P=\deg P' =\deg P''=1$. In other words, the condition that
at most two  of the three lines determined by $p$, $p'$ and $p''$ are in $\A$ is necessary for our proof of Theorem \ref{thm3} to hold.

\end{ex}

\section{The proofs of Theorems \ref{thm1G}, \ref{thmG} and of Corollary \ref{cor1G}} 

\subsection{Proof of Theorem \ref{thm1G}} 
For the morphism $u_0$ one may refer to \cite[Theorem 5.1]{DIS},
and the fact that $v_0$ is well defined can be checked by a direct computation. The claim that any $\delta \in D_0(f)$ has the form $\delta=(xa_1,b,c)$ with $a_1 \in S$ follows from the observation
that $f=xf'$ with $f'$ and $f_x=f'+xf_x'$ not divisible by $x$, while both $f_y=xf'_y$ and $f_z=xf'_z$ are divisible by $x$.
To prove that the sequence \eqref{eq1} is exact with these choices for $u_0$ and $v_0$, we proceed as follows.
Note that $\delta=(xa_1,b,c)$ is in $\ker v_0$ if and only if
$$(d-1)\overline b+\overline a_1y=0 \text{ and } (d-1)\overline c+\overline a_1z=0.$$
It follows that there is a polynomial $b_2 \in \C[y,z]$ such that
$$\overline b=b_2y, \ \overline c=b_2z \text{ and } \overline a_1=-(d-1)b_2.$$
Hence we can write
\begin{equation}
\label{eq3}
a=x(xa_2-(d-1)b_2), \  b=xb_1+yb_2 \text{ and } c=xc_1+zb_2
\end{equation}
for some polynomials $a_2,b_1,c_1 \in S$.
With this notation, the condition $\delta(f)=0$ implies that
$\delta'(f')=0$, where
$$\delta'=\left( \frac{da_2x}{d-1}-db_2, \frac{a_2y}{d-1}+b_1,  \frac{a_2z}{d-1}+c_1\right).$$
A direct computation shows that $u_0(\delta')=\delta$, and hence
$$\ker v_0 \subset \im u_0.$$
The other facts to check for the sequence \eqref{eq1} to be exact are obvious, hence this completes the proof of Theorem \ref{thm1}.

\subsection{Proof  of Corollary \ref{cor1G}} 
Note that the module $D_0(f'')$ is generated  by a single syzygy $\delta''$ in degree $k=r-1$. Indeed, it is enough to take $\delta''$ the syzygy obtained from the Koszul syzygy
 $$f''_zf''_y-f''_yf_z''=0$$
 which has degree $d-1$ by factoring out the common divisor of $f'_y$ and $f'_z$ which has degree
 $d-k$. 
 
 It follows that all the generators of $D_0(f)$ of degree $<r-1$ are mapped to 0 by $v_0$, and hence they are in $\im u_0$. This remark proves the claim $(1)$.
 If we assume that $d_s<k-1$, then it follows that the morphism of graded $S$-modules $v_0$ is trivial and hence we get an isomorphism $D_0(f)=D_0(f')(-1)$.
In particular, for any integer $k$, we have an isomorphism 
$$D_0(f')_{k-1}=D_0(f)_k$$
of $\C$-vector spaces.
To show that this is impossible, we use the exact sequence in \cite[Corollary 2.5]{POG}. With the notation there, the divisor $D'$ on $C_2=L$ is a point, and hence for $k$ large, one has 
$$H^0(C_2, \OO_{C_2}(D+(k-1)D') \ne 0.$$
On the other hand, for $k$ large, one has for $f'=f_1$ the following vanishing
$$N(f_1)_{k-e_2+e_1-1}=0,$$
since the Jacobian module $N(f_1)$ is Artinian, see \cite{Se}. This contradiction proves our claim $(2)$.

\subsection{Proof of Theorem \ref{thmG}} 

From now on in this section we assume that $\A=\A' \cup L$, where $\A':f'=0$ is line arrangement obtained by deleting the line $L$ from $\A$. 
Let $\A'':f''=0$ be the multiarrangement on $L$ defined by the restriction of $f$ to $L$.
We choose the coordinates such that $L:x=0$ and, then $f=xf'$ and $f''=\overline f'$,
In order to use Theorem \ref{thm1G}, we first  see how the local derivations $\tilde D_p$ are mapped under the morphisms $u_0$ and
$v_0$.

\begin{prop}
\label{propG1}
Consider an arrangement $\A$ of $d$ lines in $\PP^2$ such that $L:x=0$ is a line in $\A$ and set $\A'=\A \setminus \{L\}$. For a multiple point
$p \in \A'$ we denote by $\tilde D_p' \in D_0(\A')$ the local derivation associated to $p$, and by $\tilde D_p \in D_0(\A)$ the local derivation associated to $p \in \A$. Then the following two cases are possible.

\begin{enumerate} 

\item If $p \in L$, then $u_0(\tilde D_p')=x\tilde D_p$.

\item If $p \notin L$, then $u_0(\tilde D_p')=\tilde D_p$.

\end{enumerate} 

\end{prop}

\proof
To prove the claim $(1)$, we may assume that $p=(0:1:0)$ and $f=xf'$.
Then, with obvious notation, we have $f_1=xf_1'$ and $f_2=f_2'$.
It follows that
$$\tilde D_p'=f_2\partial_y-\frac{f_{2,y}}{d-1}E$$
and hence
$$u_0(\tilde D_p')=    x\tilde D_p'+\frac{xf_{2,y}}{d(d-1)}E=x\tilde D_p.$$
To prove the claim $(2)$, we may assume that $p=(1:0:0)$ and $f=xf'$.
In this case we have $f_1=f_1'$ and $f_2=xf_2'$.
It follows that
$$\tilde D_p'=f_2'\partial_x-\frac{f_{2,x}'}{d-1}E$$
and hence
$$u_0(\tilde D_p')=  x\tilde D_p'-\frac{(d-1)f_2'-xf_{2,x}}{d(d-1)}E=\tilde D_p.$$
\endproof

\begin{prop}
\label{propG2}
Consider an arrangement $\A$ of $d$ lines in $\PP^2$ such that $L:x=0$ is a line in $\A$ and set $\A'':f''=0$ to be the multiarrangement obtained by restriction of $\A$ to $L$.
For a multiple point
$p \in \A$ we denote by $\tilde D_p \in D_0(\A)$ the local derivation associated to $p$. Then the following two cases are possible.

\begin{enumerate} 

\item If $p=(0:-\gamma_p: \be_p) \in L$, then we set $f=f_1f_2$ and $f''=\overline f_2=f_2(0,y,z)$ as above. With this notation, one has 
$$v_0(\tilde D_p)=\frac{1}{d-1} (yA_p-(d-1)\gamma_p\overline f_2, zA_p+(d-1)\be_p\overline f_2),$$
where $A_p=\gamma_p \overline f_{2,y}-\be_p\overline f_{2,z}$.

\item If $p \notin L$, then $v_0(\tilde D_p)=0$.

\end{enumerate} 

\end{prop}

\proof
To prove the first claim, notice that
$$\tilde D_p=f_2(0,-\gamma_p,\be_p)+\frac{A_p}{d}(x,y,z).$$
The formula for $v_0(\tilde D_p)$ follows now from the definition of $v_0$ in Theorem \ref{thm1G}. The second claim follows from the exact sequence in Theorem \ref{thm1G} and Proposition \ref{propG1}, claim $(2)$.
\endproof

\subsection{The proof of Theorem \ref{thmG}} 

With the notation from Propositions \ref{propG2}, let $K=L \cap \A'$ and let $k=|K|$ be the number of points in this intersection. Since $\A$ is not a pencil, we get $k\geq 2$. Let $p_1, \ldots, p_k$ be the points in $K$ and assume that
$$p_i=(0:-\gamma_i:\be_i) \text{ for } i=1, \ldots , k.$$
Then one has $f=xf'$ and
\begin{equation}
\label{eqG1}
f''(y,z)=f'(0,y,z)=\prod_{i=1}^k\ell_i^{m_i-1},
\end{equation}
where $\ell_i=\be_i y+\gamma_i z$ and $m_i=m_{p_i}$. In particular, we have
\begin{equation}
\label{eqG2}
\sum_{i=1}^k(m_i-1)=d-1.
\end{equation}
It follows that the $S$-module $D_0(f'')$ has one generator $\delta''$, namely the one obtained from the Koszul derivation
$$\delta _K=f''_z\partial_y-f"_y \partial_z$$
by removing the common factor
$$G.C.D.(f''_y,f''_z)=\prod_{i=1}^k\ell_i^{m_i-2}.$$
It follows that
\begin{equation}
\label{eqG3}
\delta''=\sum_{i=1}^k(m_i-1)(\prod_{s \ne i}\ell_s)(\gamma_i \partial_y-\be_i\partial_z)=B''\partial_y +C''\partial_z.
\end{equation}
Now we consider the derivation $\tilde D_{p_1}$ and find out the relation between $v_0(\tilde D_{p_1})$ and $\delta''$.
To use Proposition \ref{propG2} $(1)$, we first notice that all the lines in $\A$ which do not pass through $p_1$, intersect $L$ in the points $p_s$ with $s \ne 1$. It follows that, if $f_2=0$ is the equation of the union of these lines in $\A$, one has
\begin{equation}
\label{eqG4}
g_1=f_2(0,y,z)=\prod_{s=1,\ldots,k; s \ne 1}\ell_s^{m_s-2},
\end{equation}
where the index $1$ in $g_1$ says that this polynomial is constructed with respect to the point $p_1 \in K$.
Proposition \ref{propG2} $(1)$ yields the following equality, where the constant factor $1/(d-1)$ was omitted since it plays no role
$$v_0(\tilde D_{p_1})=g_1(B_1\partial_y+C_1\partial_z)$$
where
$$B_1=\gamma_1 y\sum_{ s \ne 1}\left((m_s-1)\be_s\prod_{t \ne 1, t \ne s}\ell_t\right) -  \be_1y \sum_{ s \ne 1}\left((m_s-1)\gamma_s\prod_{t \ne 1, t \ne s}\ell_t\right)-(d-1)\gamma_1\prod_{ s \ne 1}\ell_s$$
and
$$C_1=\gamma_1 z\sum_{ s \ne 1}\left((m_s-1)\be_s\prod_{t \ne 1, t \ne s}\ell_t\right)-  \be_1z \sum_{ s \ne 1}\left((m_s-1)\gamma_s\prod_{t \ne 1, t \ne s}\ell_t\right)+(d-1)\be_1\prod_{ s \ne 1}\ell_s,$$
where $s,t$ take all values from $1$ to $k$ except those indicated in the sums and respectively in the products above. Using \eqref{eqG2} we can rewrite $B_1$ as follows
$$B_1=\gamma_1 y\sum_{ s \ne 1}\left((m_s-1)\be_s\prod_{t \ne 1, t \ne s}\ell_t\right) -  \be_1y \sum_{ s \ne 1}\left((m_s-1)\gamma_s\prod_{t \ne 1, t \ne s}\ell_t\right)-$$
$$-\gamma_1 \sum_{ s \ne 1}\left((m_s-1)\ell_s \prod_{t \ne 1, t \ne s}\ell_t\right)-\gamma_1(m_1-1)\prod_{s \ne 1}\ell_s.$$
Note that
$$\gamma_1\be_sy-\gamma_1\ell_s=-\gamma_1\gamma_sz.$$
Hence after cancelation of terms in the first and third sum above we get
$$B_1=-\ell_1 \sum_{ s \ne 1}\left((m_s-1)\gamma_s \prod_{t \ne 1, t \ne s}\ell_t\right)- \gamma_1(m_1-1)\prod_{s \ne 1}\ell_s.$$
It follows that $B_1=-B''$, and a similar computation yields $C_1=-C''$.
In conclusion, one has
\begin{equation}
\label{eqG5}
 v_0(\tilde D_{p_1})=g_1(B_1\partial_y+C_1\partial_z)=-g_1\delta''.
\end{equation}
By symmetry, if we set for any $i \in [2,k]$,
\begin{equation}
\label{eqG6}
g_i=\prod_{s=1,\ldots,k; s \ne i}\ell_s^{m_s-2},
\end{equation}
we get the equality
\begin{equation}
\label{eqG7}
 v_0(\tilde D_{p_i})=-g_i\delta'',
\end{equation}
where as above the constant factor $1/(d-1)$ is omitted.
To continue the proof of Theorem \ref{thmG}, we show an inclusion of ideals in the ring $\C[y,z]$, namely
\begin{equation}
\label{eqG8}
 (y,z)^{d-2-k} \subset I=(g_1, \ldots,g_k).
\end{equation}
Note that in view of \eqref{eqG2} we have.
$$\sum_{i=1}^k(m_i-2)=d-1-k$$
and that we may assume $\ell_1=y$ and $\ell_2=z$.
We prove this claim by induction on $k$. If $k=2$,  the claim is clear. Indeed, any monomial $y^sz^t $ of degree $s+t=d-4$ cannot satisfy
$$s \leq m_1-3 \text{ and } t \leq m_2-3,$$
since then $s+t \leq (m_1-3)+(m_2-3)=d-5$, a contradiction.
Hence $y^sz^t$ is divisible either by $g_2$ or by $g_1$, and hence
$y^sz^t \in I$ in this case.
Assume now the claim \eqref{eqG8} proved for any set of $k-1$ linear forms.
By this hypothesis, it follows that we have the inclusions of ideals
$$y^{m_1-2}(y,z)^{\sum_{i \ne 1}(m_i-2)-1} \subset I  \text{ and } z^{m_2-2}(y,z)^{\sum_{i \ne 2}(m_i-2)-1} \subset I.$$

If we have a monomial $y^sz^t$ such that $s+t=d-2-k$, then two cases are possible. Either $s \geq m_1-2$ and then 
$$y^sz^t \in y^{m_1-2}(y,z)^{\sum_{i \ne 1}(m_i-2)-1},$$
or $t \geq m_2-2$ and then
$$y^sz^t \in z^{m_2-2}(y,z)^{\sum_{i \ne 2}(m_i-2)-1}.$$
Indeed, if $s \leq  m_1-3$ and $t \leq m_2-3$ we get
$$\sum_{i=1}^k(m_i-2)=d-1-k=s+t+1\leq (m_1-2)+(m_2-2)-1.$$
This is a contradiction, since $m_i -2 \geq 0$ for any $i$.
Hence the proof of the inclusion \eqref{eqG8} is completed.

Now we prove Theorem \ref{thmG} by induction on $d=|\A| \geq 3$.
Since we assume that $\A$ is not a pencil, the starting point of the induction is $\A:f=xyz=0$. Then it is known that $\A$ is free with exponents $(1,1)$. In particular $\dim D_0(\A)_1=2$ and any basis $\rho_1,\rho_2$ of the vector space $D_0(\A)_1$ freely generate the $S$-module $D_0(\A)$.
For the double point $p=(1:0:0)$ (resp. $p'=(0:1:0)$) we get
$$\tilde D_p=x\partial_x-\frac{1}{3}E \in D_0(\A)_1$$
and respectively
$$\tilde D_{p'}=y\partial_y-\frac{1}{3}E \in D_0(\A)_1.$$
These two local derivations are clearly independent, form a basis of the vector space $D_0(\A)_1$, and hence generate $D_0(\A)_k$ for any $k \geq d-2=1$.
Assume the claim in Theorem \ref{thmG} holds for all arrangements $\A'$ with
$|\A'|<d$ and consider the exact sequence \eqref{eq2G2}, which yields
\begin{equation}
\label{eq2G}
0 \to D_0(f')_{d-3} \stackrel{u_0} \longrightarrow D_0(f)_{d-2} \stackrel{v_0} \longrightarrow D_0(f'')_{d-2}.
\end{equation}
By the induction hypothesis, $D_0(f')_{d-3}$ is generated by local derivations corresponding to the multiple points of $\A'$. Proposition \ref{propG1} says that $\im u_0$ is contained in the graded submodule $M \subset D_0(\A)$ generated by the local derivations corresponding to the multiple points of $\A$. 
On the other hand, the equalities \eqref{eqG7} and the inclusion \eqref{eqG8} imply that the morphism
$$M_{d-2} \stackrel{v_0} \longrightarrow D_0(f'')_{d-2}=S_{d-k-1}\rho''$$
obtained by restriction of $v_0$ to $M_{d-2}$ is surjective. The exact sequence \eqref{eq2G} shows that $M_{d-2}=D_0(f)_{d-2}$, which is enough to prove
Theorem \ref{thmG} in view of our Remark \ref{rkG} above.

\section{Local derivations, Bourbaki ideals and the proof of Theorem \ref{thm4} and Proposition \ref{prop40}} 

We recall first the construction of the Bourbaki ideal $B(C,\rho_1)$ associated to a degree $d$ reduced curve $C:f=0$ and to a minimal degree non-zero syzygy $\rho \in D_0(f)$, see \cite{DStJump}
as well as \cite{JNS} for a recent, more complete approach.
For any choice of the syzygy $\rho=(a,b,c)$ with minimal degree $d_1$, we have a morphism of graded $S$-modules
\begin{equation} \label{B1}
S(-d_1)  \xrightarrow{u} D_0(f), \  u(h)= h \cdot \rho.
\end{equation}
For any homogeneous syzygy $\rho'=(a',b',c') \in D_0(f)_m$, consider the determinant $\Delta(\rho')=\det M(\rho')$ of the $3 \times 3$ matrix $M(\rho')$ which has as first row $x,y,z$, as second row $a,b,c$ and as third row $a',b',c'$. Then it turns out that $\Delta(\rho')$ is divisible by $f$, see \cite{Dmax}, and we define thus a new morphism of graded $S$-modules
\begin{equation} \label{B2}
 D_0(f)  \xrightarrow{v}  S(d_1-d+1)   , \  v(\rho')= \Delta(\rho')/f,
\end{equation}
and a homogeneous ideal $B(C,\rho) \subset S$ such that $\im v=B(C,\rho)(d_1-d+1)$.
We have  
\begin{equation} \label{B3}
 B(C,\rho)=S \text{ if and only if } C \text{ is a free curve}.
\end{equation}
It is known that the ideal $B(C,\rho)$, when $C$ is not a free curve, defines a $0$-dimensional subscheme in $\PP^2$, which is locally a complete intersection, see  for details \cite[Theorem 5.1]{DStJump}.

Let $\A:f=0$ be a line arrangement and fix a minimal degree non-zero syzygy $\rho \in D_0(f)$, say $\rho=(a,b,c)$ and $\deg \rho=d_1$.
For any multiple point $p \in \A$ we consider the corresponding local derivation $\tilde D_p$ and define a polynomial $g_p \in S$, such that one has
$$g_p=v(\tilde D_p) \in B(\A,\rho).$$
Using the formula for $\tilde D_p$ given in Theorem \ref{thm1} we see that
$$\Delta(\tilde D_p)=f_{2p}\Delta_p$$
where $\Delta_p$ is the determinant  of the $3 \times 3$ matrix $M_p$ which has as first row $x,y,z$, as second row $a,b,c$ and as third row $u,v,w$ if $p=(u:v:w)$. It follows that
\begin{equation} \label{B4}
 g_p=\frac {\Delta(\tilde D_p)}{f}=\frac{\Delta_p}{f_{1p}}.
\end{equation}
One has the following result.
\begin{prop}
\label{prop4}
With the above notation, for each multiple point $p \in \A$ of multiplicity $m_p$ satisfying
 \eqref{E1}
the corresponding polynomial $g_p \in S$ is homogeneous of degree
$d_1+1-m_p \geq 0$ when $\tilde D_p$ is not a constant multiple of $\rho$. In particular $g_p \ne 0$ in this situation.
\end{prop}

\proof
Note that $g_p=0$ is equivalent to $\Delta(\tilde D_p)=0$. This last condition is equivalent to the fact that the rows in the matrix $M(\Delta(\tilde D_p))$ are linearly dependent over the field of franctions of the polynomial ring $S$. This occurs if and only if there is a nontrivial relation
$$P E +Q D_{\rho} +R \tilde D_p=0,$$
where $P,Q,R \in S$ and $E$ is as above the Euler derivation.
If we apply this relation to $f$, we get $d P f=0$, and hence $P=0$.
Suppose $\tilde D_p$ is not a constant multiple of the minimal degree derivation $\rho$. Since $\tilde D_p$ is a primitive derivation, it follows that $\tilde D_p$ is not a polynomial multiple of $\rho$. Then the result follows, since the quotient
$D_0(f)/\im u$, with $u$ as in \eqref{B1}, is a torsion free $S$-module, being isomorphic to $\im v$. Hence $Q=R=0$, a contradiction which implies that $g_p \ne 0$.
If $g_p \ne 0$, then the claim $\deg g_p=d_1+1-m_p$ follows from \eqref{B4}. 
\endproof

\subsection{Proof of Theorem 
\ref{thm4}}

Let $q \in \A$ be a point of multiplicity $m$. If $m>d/2$, then the associated derivation $\tilde D_q$ is a derivation in $D_0(f)$ of minimal degree $r=d-m \leq (d-1)/2$, see \cite[Theorem 1.2]{Mich}. The claim about freeness of $\A$ follows from \cite{dPCTC,Dmax}. This proves the claim $(1)$. To prove the claim $(2)$, we apply again  \cite[Theorem 1.2]{Mich} and obtain the claim for the case
$d_1=m-1$, as well as the inequalities $m \leq d_1 \leq d-m$.

The proof for $(3)$ is more subtle.
The case $m=2$ is clear, since we have $d \geq 4$ and such a nodal arrangement is never free. On the other hand $I(\rho)=0$ in this case.
Assume from now on that $m \geq 3$.  Consider first the case when $\A$ is free, with exponents $d_1 \leq d_2$. Then $d_1 \geq 2$  by $(1)$ and $(2)$ and hence
$$d_2=d-1-d_1\leq d-3.$$
It follows that the second generator $\rho_2$ of $D_0(f)$ has degree
$d_2 \leq d-3$ and hence we have the equality
$$S_{d-3-d_1}\rho_1 + S_{d-3-d_2}\rho_2 = D_0(f)_{d-3}.$$
On the other hand, since $\A$ is free, we know that $v(\rho_2)=\lambda$
is a non zero constant. This implies that
$$v(D_0(f)_{d-3})=S_{d-3-d_2} \subset B(\A,\rho).$$
Any element $D' \in D_0(f)_{d-3}$ can be written by Theorem \ref{thm2} as a sum
$$ D'=\sum A_p \tilde D_p,$$
where $A_p \in S_{m_p-3}$ and $p$ runs through all the multiple points of $\A$ satisfying $m_p \geq 3$. Hence
$$v(D')=\sum A_p v(\tilde D_p)=\sum A_p g_p.$$
Since in this sum
$$0\leq \deg g_p=d_1+1-m_p \leq \frac{d+1}{2}-m_p,$$
it follows that the points $p$ in the above sum satisfy \eqref{E1}.
This shows that 
$$S_{d-3-d_2} \subset I(\rho),$$
and hence the zero set of the ideal $I(\rho)$ is empty in $\PP^2$.
Conversely, if the zero set of the ideal $I(\rho)$ is empty in $\PP^2$,
then the same holds for the ideal $B(\A,\rho)$, since by definition
$$I(\rho) \subset B(\A,\rho).$$
But the ideal $B(\A,\rho)$ being saturated, see \cite[Theorem 5.1]{DStJump}, this is possible if and only if
$B(\A,\rho)=S$. Using \eqref{B3}, we conclude that $\A$ is free, and this completes the proof of Theorem \ref{thm4}.

\subsection{Proof of Proposition
\ref{prop40}}
We choose $\rho=\tilde D_q$ as the minimal degree derivation and
the coordinates such that $q=(1:0:0)$ and $p=(0:0:1)$.
Then $\Delta(\tilde D_p)$ is equal to the determinant of the $3 \times 3$ matrix having as first row $x,y,z$, as second row $(f_{2q},0,0)$ and the third row $(0,0,f_{2p})$, and hence, up to a sign, one has
$$\Delta(\tilde D_p)=yf_{2p}f_{2q}.$$
\noindent{\bf Case 1:} The line $L(p,q):y=0$ is in $\A$. Then one can write $f_{1p}=yf_{1p}'$, such that all the lines defined by $f_{1p}'=0$ are in $\A$ and do not contained $q$. It follows that 
$$f_{2q}=f_{1p}'h,$$
where $h=0$ is the union of lines containing neither $p$ nor $q$. It follows that $h=0$ is an equation for the arrangement $\B_p$ defined above, and hence $h=h_p$ up to a constant factor.
On the other hand we have
$$g_p= \frac{\Delta(\tilde D_p)}{f}=\frac{yf_{2p}f_{2q}}{f_{1p}f_{2p}}=
\frac{yf_{2p}f_{1p}'h}{f_{1p}f_{2p}}=h,$$
which proves our claim in this case.

\medskip

\noindent{\bf Case 2:} The line $L(p,q):y=0$ is not in $\A$.
Then one can write as above $f_{2q}=f_{1p}h$, where
where $h=0$ is the union of lines containing neither $p$ nor $q$.
 It follows that $yh=0$ is an equation for the arrangement $\B_p$ defined above, and hence $h_p=yh$ up to a constant factor.
On the other hand we have
$$g_p= \frac{\Delta(\tilde D_p)}{f}=\frac{yf_{2p}f_{2q}}{f_{1p}f_{2p}}=
\frac{yf_{2p}f_{1p}h}{f_{1p}f_{2p}}=yh,$$
which proves our claim in this case and completes thus the proof of Proposition \ref{prop40}.

\begin{rk}
\label{rk4}
Note that if $m \leq d/2$ and
$$\tau(\A)=(d-1)^2-m(d-m-1),$$
then the arrangement $\A$ is free with exponents $(m,d-m-1)$.
Indeed, using  \cite{dPCTC,Dmax} it follows that $d_1 \leq m$.
Theorem \ref{thm4} $(2)$, implies that $d_1 \geq m-1$, and the equality
$d_1=m-1$ is not possible. Indeed, again by Theorem \ref{thm4} $(2)$,
$d_1=m-1$ implies that
$$\tau(\A)=(d-1)^2-(m-1)(d-m)>(d-1)^2-m(d-m-1),$$
a contradiction. Therefore, the simplest case where the freeness cannot be decided easily by the combinatorics seems to be the case
$$d_1=m+1 \text{ and } 2m+3 \leq d.$$
For instance, the free arrangement $\A$ with exponents $(m+1,m+1)$ and the nearly free arrangement $\B$ with exponents $(m,m+3,m+3)$
satisfy $\tau(\A)=\tau(\B)$. We do not know whether such a pair of arrangements $\A,\B$ exists, such that $\A$ and $\B$ have the same combinatorics.
\end{rk}

\section{Proof of Theorems 
\ref{thm1000} and \ref{thm100}  and of Corollaries \ref{cor101} and \ref{cor1000}}
We start with the following preliminary results.
\begin{prop}
\label{propthm10}
Let $\A:f=0$ be a line arrangement and $\rho=(a,b,c) \in D_0(f)$.
Then the following holds.

\begin{enumerate} 

\item If $L:\ell=\al x +\be y + \gamma z=0$ is one of the lines in $\A$, then
$$\sigma=\al a + \be b + \gamma c$$
is divisible by $\ell$.
\item If $p \in \A$ is a point of multiplicity at least $2$, then 
the vector $\rho(p)=(a(p),b(p),c(p))$ is either zero, or a multiple of $p$.

\item Let $q \in \A$ be a point of multiplicity at least $2$ and such that
$p \ne q$ and the line determined by $p$ and $q$ is not in $\A$.
Then $g_p(q)=0$.

\end{enumerate} 
\end{prop}
Note that the claims $(1)$ and $(2)$ are well known, see for instance
\cite[Proposition 3.6]{To} and its proof.

\proof

To prove $(1)$, we note that if we set $f= \ell f'$, then 
$D_{\rho}(f)$ is a sum of terms divisible by $\ell$ plus the sum
$f'\sigma$. It follows that $\sigma$ is divisible by $\ell$.
If we denote by $L^0$ the hyperplane in $\C^3$ given by $\ell=0$, then
the claim $(1)$ says exactly that $D_{\rho}$, regarded as a vector field on $\C^3$, is everywhere tangent to $L^0$.

To prove $(2)$, note that if $p \in L_1 \cap L_2$, then the tangent vector
$D_{\rho}(p)$ is in the line $L_1^0 \cap L_2^0=\C \cdot p$.
To prove $(3)$, we recall that 
$$\Delta_p=g_p f_{1p},$$
see \eqref{B4}. If we evaluate this equality at the point $q$ we have
$f_{1p}(q) \ne 0$, since the line determined by $p$ and $q$ is not in $\A$. Moreover we have $\Delta_p(q)=0$. Indeed, the first two rows in the matrix used to compute $\Delta_p(q)$ are
$v_1=(u,v,w)$ and $v_2=(a(q),b(q),c(q))$, where we assume that $q=(u:v:w)$.
The vector $v_2$ is a multiple of the vector $v_1$ by the claim $(2)$,
and hence this determinant is zero.
\endproof

\begin{lem}
\label{lem10}
Let $\A:f=0$ be an arrangement of $d$ lines in $\PP^2$ and $p \in \A$  a point of multiplicity $m= m(\A)$, the maximal multiplicity of points in $\A$. The projection $\phi_p: \PP^2 \setminus \PP^1$ induces a locally trivial fibration $\psi_p: M(\B) \to B$ with fiber $F=\PP^1 \setminus Q_F$, where $Q_F$ is a finite set consisting of $d-m+1$ points,
base $B=\PP^1 \setminus Q_B$, where $Q_B$ is a finite set consisting of $b \geq m$ points, and $M(\B)$ is the complement in $\PP^2$ of the line arrangement obtained from $\A$ by adding all the lines $L_q \notin \A$, which connect $p$ to a multiple point $q \ne p$ in $\A$.
Then the following holds.

\begin{enumerate} 

\item Let $E$ be the set of lines added to $\A$ to get the new arrangement $\B$ and let $e=|E|$. Then $b=m+e$.

\item For each line $L_q \notin \A$, we denote by $n_q$ the number of points in the intersection $I_q=L_q \cap \A$. Then the Euler number of the complement $M(\A)$ of the line arrangement $\A$ in $\PP^2$ is given by
$$E(M(\A))=E(M(\B))+\sum_qE(L_q \setminus I_q)=
(m+e-2)(d-m-1)+\sum_q(2-n_q).$$

\end{enumerate} 
\end{lem}
\proof
The fact that $\psi$ is a locally trivial fibration (in the strong complex topology) is clear, since the set $E$ of lines we add to $\A$ consists exactly of the lines $L \notin \A$ such that $|L \cap \A|$ has a strictly smaller value than the generic one, which is $d-m+1$. This remark proves also the claim (1). To prove the claim (2), we just recall that the Euler number is additive with respect to constructible partitions.
\endproof

\subsection{Proof of Theorem 
\ref{thm1000}}
The claim (1) and the claims before it follow directly from Proposition \ref{propthm10}.
Note that the inequality $d_1<d-m_p$ implies that $\tilde D_p$ is not a constant multiple of the minimal degree derivation $D_{\rho}$, and hence $g_p \ne 0$.

The first claim in (2)  also follows directly from Proposition \ref{propthm10}. To prove the rest of the claim (2)(a) we have to consider three cases. To simplify the notation we set $m=m_p$.

\medskip

\noindent {\bf Case 1} ($L\notin \A$ and $p \notin L$). Then for any line $L_q$ to be added to $\A$ as in Lemma \ref{lem10}, the point $q$ is in $L$ and the line $L_q$ contains no other multiple points of $\A$, except $p$ and $q$.
It follows that, in the notation from Lemma \ref{lem10}, we have
$$n_q=d-m-m_q+2.$$
It follows that
$$E(M(\A))=(m+e-2)(d-m-1)+\sum_q(m+m_q-d)=(m-2)(d-m-1)+\sum_q(m_q-1).$$
On the other hand, using again that the Euler number is additive with respect to constructible partitions, we get
$$E(M(\A))=E(\PP^2)-E(\A)=3-(2-(d-1)(d-2)+\tau(\A)),$$
see for instance \cite[Equation (4.5)]{DHA}.
By comparing the last two equalities, we get
$$\tau(\A)=(d-1)^2-m(d-m-1)-\left( \sum_q(m_q-1) -(d-2m)\right).$$
Using Remark \ref{rk4} and \cite{dPCTC}, since $d_1=m$, it follows that
$$  \sum_q(m_q-1) -(d-2m) \geq 0$$
and the equality holds if and only if $\A$ is free.
On the other hand, it is clear that
$$ \sum_q(m_q-1)= \sum_q m_q-e \leq \sum_{q'}m_{q'}-e=d-m-e,$$
where the second sum is over all points $q' \in L \cap\A_2$  and $m_{q'}'$ denotes the multiplicity of $q'$ in $\A_2$.
Note that the points $q$ in the above sum are exactly the multiple points of the arrangement $\A_2:f_{2p}=0$ which are not connected to $p$ by a line in $\A$, and such a point $q$ has the same multiplicity in $\A$ and in $\A_2$.
Bezout's Theorem implies that
$$\sum_{q'}m_{q'}'=d-m$$
where the sum is over all points $q' \in L \cap\A_2$.
In conclusion, in this case we have
\begin{equation} \label{I2}
d-2m \leq \sum_q(m_q-1) \leq d-m-e.
\end{equation} 
This shows in particular that one has
\begin{equation} \label{I3}
 1\leq e \leq m.
\end{equation} 

\medskip

\noindent {\bf Case 2} ($L \notin \A$ and $p \in L$). Let $n_L=|L \cap \A_2|$ as above and hence $|I|=n_L+1$, where $I=L \cap \A$. With the notation from Lemma \ref{lem10} we have $e=1$,
$E(L\setminus I)=1-n_L$ and hence
$$E(\A)=(m-1)(d-m-1)+1-n_L.$$
As above we conclude that
$$\tau(\A)=(d-1)^2-m(d-m-1)-(m-n_L).$$
Using Remark \ref{rk4} and \cite{dPCTC}, since $d_1=m$, it follows that
\begin{equation} \label{I1}
 m -n_L \geq 0
\end{equation} 
 and the equality holds if and only if $\A$ is free.
 Note that the equality
$$d=\sum_{q \in I}m_q$$
can be rewritten as
$$d=\sum_{q \in I}(m_q-1)+n_L+1=m-1+\sum_{q \in I, q \ne p}(m_q-1)+n_L+1.$$
It follows that
$$n_L = d-m-\sum_{q\in I, q \ne p}(m_q-1).$$
Hence, the inequality \eqref{I1} becomes
$$\sum_{q\in I, q \ne p}(m_q-1) \geq d-2m,$$
hence the same as in Case 1.

\noindent {\bf Case 3} ($L\in \A$).

Now the line $L$ is in $\A$,  and it is clear that $L$ is in fact a line in $\A_2$ as well, since otherwise $\A$ would be semisimple with $p$ a modular point. This would imply $d_1=m-1$, a contradiction. Hence in this case $p \notin L$ and for any point $q \in L \cap \A_2$ with multiplicity in $\A$ satisfying $m_q \geq 2$, if the line $L_q$ is not in $\A$, then $L_q$ has only $p$ and $q$ as multiple points. It follows that one has, 
$$|L_q \cap \A|=d-m-m_q+2,$$
exactly as in Case 1 above.
Using the same approach as in Case 1, we see that
$$\tau(\A)=(d-1)^2-m(d-m-1)-\left( \sum_q(m_q-1) -(d-2m)\right).$$
Since $d_1=m$,  Remark \ref{rk4} and \cite{dPCTC} imply  that
$$  \sum_q(m_q-1) -(d-2m) \geq 0$$
and the equality holds if and only if $\A$ is free.
This 3 cases may occur for one and the same line arrangement $\A$, see Example \ref{ex100} below, corresponding to various choices of the line $L$.

Suppose now that 
$$  \sum_q(m_q-1) >d-m. $$
Then the arrangement $\A$ is not free, and hence its second exponent satisfies $d_2 \geq d-m$. On the other hand, since $\tilde D_p$ is a primitive derivation of degree $d-m > d_1=m$, it follows that
$d_2\leq d-m$. Therefore, $d_2=d-m$ and $\A$ is a plus-one generated arrangement, see \cite[Theorem 2.3]{3syz}, with exponents
$d_1=m$, $d_2=d-m$ and $d_3$. To determine $d_3$ we use the formula
$$\tau(\A)=(d-1)^2-d_1(d-1-d_1)-(d_3-d_2+1),$$
which comes from \cite[Proposition 2.1 (4)]{3syz} using the relation
$d_1+d_2=d$.

\endproof

\begin{rk}
\label{rkex10}
 If the line $L$ in Theorem \ref{thm1000} is not in $\A$, then
$$\sum_q m_q \leq d-m_p,$$
where $m_q$ denote the multiplicity of point $q$  in $\A$. This is just a reformulation of the second inequality in \eqref{I2} in Case 1 in the above proof. If this line $L$ is in $\A$, then one has the (usually) weaker inequality
$$\sum_q (m_q-1) \leq d-m_p-1.$$
These claims are  direct consequences of Bezout's Theorem, and they hold even when $m_p <m(\A)$, see Example \ref{ex10}, (ii) and (iii).
\end{rk}

\begin{ex}
\label{ex10}
(i) Consider the Hessian line arrangement
$$\HH: xyz\left( (x^3+y^3+z^3)^3-27x^3y^3z^3\right )=0$$
which is free and it has $d_1=4=m(\A)$, see \cite[Example (8.6) (iii)]{DHA}. Note that $\HH$ is the union of 4 line arrangements $\A_1,\A_2, \A_3$ and $\A_4$, each one consisting of one triangle. More precisely, 
$$\A_j: x^3+y^3+z^3-3\epsilon ^jxyz=0,$$
for $j=1,2,3$ and $\epsilon=\exp(2\pi i/3)$, and $\A_4:xyz=0$.
The point $p=(0:1:-1)$, has multiplicity $m_p=4=m(\A)$, and the 4 lines passing through $p$ are $x=0$, $x+y+z=0$, $\epsilon x+y+z=0$
and $\epsilon^2 x+y+z=0$, hence exactly one line in each triangle $\A_j$. It follows that in each triangle $\A_j$ there is a double point of $\HH$ not situated on these 4 lines. These points are given by
$$(1:\epsilon: \epsilon), (1:\epsilon^2: \epsilon^2), (1:1:1) \text{ and } (1:0:0).$$
They are situated on the line $L:y-z=0$ which is not in $\HH$. This line is as in Case 1 in the above proof, and one has
$$\sum_{q\in Q}(m_q-1)=4=d-2m_p$$
as predicted by Theorem \ref{thm1000}.
Next we show that there are no 5 lines  in $\HH$, containing all the multiple points of $\HH$. Indeed, if such a set of lines exists, there is a triangle $\A_k$ containing only one line $L'$ in this set. It follows that the double point in this triangle $\A_k$ which is not on the line $L'$, is not contained in the union of the 5 lines.

\noindent (ii) Consider the full monomial arrangement
$$\F_3:xyz(x-y)(y-z)(x-z)=0.$$
Then $d_1=2$ and let $p$ be the point $(0:1:1)$ of multiplicity $2<3=m(\A)$.
The points not connected by lines in $\F_3$ to the point $p$ are the points $(1:0:1)$ and $(1:1:0)$ which are situated on the line $$L: x-y-z=0,$$
which is not a line in $\F_3$. Note that in this case
$$\sum_q m_q = 2+2=6-2= d-m_p,$$
and hence this inequality from Remark \ref{rkex10} is sharp.

\noindent (iii) Consider next the supersolvable arrangement
$$\A: xyz(x^2-y^2)(y-z)=0.$$
Then $d_1=2$ and choose $p=(-1:1:1)$, which has multiplicity $m_p=2<4=m(\A)$.
The points not connected by lines in $\A$ to the point $p$ are the points $(0:1:0)$ and $(1:1:0)$ which are situated on the line $$L: z=0,$$
which is a line in $\A$. Note that in this case
$$\sum_q (m_q-1) = 1+1 < 3=6-2-1= d-m_p-1.$$

\noindent (iv) Finally, let $\B_1$ (resp. $\B_2$) be the arrangement consisting  of $n_1\geq 2 $ lines passing through a point $a_1 \in \PP^2$ (resp. $n_2$ lines  passing through a point $a_2 \in \PP^2 \setminus \B_1$) and such that $n_1 \leq n_2$. Then the arrangement $\A=\B_1 \cup \B_2$  has  exponents $(n_1,n_2,n_1+n_2-2)$, see \cite[Proposition 1.22]{ADP}. If we choose $p=a_1$, then the only multiple point not connected by lines in $\A$ to the point $p$  is the point $a_2$. Any line through $a_2$ can be chosen as the line $L$.

\end{ex}

\subsection{Proof of Corollary
\ref{cor101}}
Choose a multiple point $p \in \A$ such that $d_1<d-m_p$, which as above implies that $g_p \ne 0$.
With the notation from the proof of Theorem \ref{thm1000}, a curve containing all the multiple points of $\A$ has equation $g_pf_{1p}=0$.
Then $\CC$ is obtained from this curve by deleting some irreducible components if they contain only multiple points contained already in the other irreducible components of $\CC$.

\subsection{Proof of Corollary
\ref{cor1000}}
If $d_1<d-m_p$, then the claim follows from Theorem \ref{thm1000} $(2)$. On the other hand, if $d_1=d-m_p$, the claim follows from
Theorem \ref{thm100} which is proved below. Note that  the inequality
$d_1>d-m_p$ is impossible, since one has $$d_1 \leq d-m \leq d-m_p,$$
where $m$ denotes the maximal multiplicity of points in $\A$, see
\cite[Theorem 1.2]{Mich}.

\subsection{Proof of Theorem 
\ref{thm100}}

If $d_1=d-m$, then we have $N \leq d_1+1=d-(m-1)$, since if we delete
$m-1$ of the lines passing through $p$, the remaining $d-(m-1)$ lines contain all the multiple points.
The equality $m-1=d_1$ implies that the $\A$ is supersolvable and such a point $p$ is modular, as shown in Theorem \ref{thm1000} (1). Hence the $m=d_1+1$ lines passing through $p$ 
contain all the multiple points.

In case $(1)$ it is known that $d_1=d-2$, see for instance \cite{Sch}
and hence $N \leq d-1$ by the claim proved above.
On the other hand $N=d-1$, since if we delete 2 lines from $\A$ their intersection point is not on the union of the remaining $d-2$ lines.

In case $(2)$ we first show that $d_1=d-3$. Using Corollary \ref{corD}, we see that $1 \leq \dim D_0(f)_{d-3} \leq 3$. Since the local derivation
$\tilde D_p$ associated with a point $p$ of multiplicity 3 is primitive, it follows that even when $\dim D_0(f)_{d-3} =3$, there is no nozero derivation of degree $<d-3$. Hence $d_1=d-3$ and hence $N \leq d-2$.
Moreover, 3 lines cannot be deleted and preserving the required property, and hence $N=d-2$.

In case $(3)$, using \cite[Theorem1.2]{Mich}, it follows that either
$$d_1 =d-m \text{
or }m-1=d_1.$$
Hence we are exactly with the hypothesis of our Theorem \ref{thm100}.

In case $(4)$ let $p$ be a modular point such that its multiplicity $m$ is the highest multiplicity of points in $\A$, see \cite[Lemma 3.1]{SS}. Then one has
$$d_1=\min \{m-1,d-m\},$$
see \cite[Proposition 3.2]{SS}.
So again we are exactly with the hypothesis of our Theorem \ref{thm100}.

\begin{rk}
\label{rk100}
(i) Even when the curve $\CC$ from Corollary \ref{cor101} is a union of lines, some of the lines may not be in $\A$, as the first arrangement in Example \ref{ex10} shows. To have another example, in the monomial arrangement
$$\A_m: f=(x^m-y^m)(y^m-z^m)(x^m-z^m)=0,$$
the minimal degree derivation $\rho$ for $m \geq 4$ is given by
$$\rho=(x(x^m-2y^m-2z^m),y(y^m-2x^m-2z^m), z(z^m-2x^m-2y^m)),$$
see \cite[Example 8.6]{DHA}.
If we take $p=(1:0:0)$, the equation for $\CC$ is the following
$$\Delta_p=yz(y^m-z^m)=0,$$
and hence $\CC$ contains two lines not in $\A_m$.

\noindent (ii) One may define another numerical invariant, namely
$N_0$ to be the minimal number of lines in $\A$, containing all the multiple points of $\A$. It is clear that
$$N_0 \geq N,$$
and Example \ref{ex10} shows that this inequality may be strict.
\end{rk}
In many cases the equality $N=d_1+1$ holds in Conjecture \ref{conj10},
but not always, as the following example shows.
\begin{ex}
\label{ex100}
 The free line arrangement $\A_k^x$ obtained by adding the line $x=0$ to the monomial arrangement $\A_k$ defined in Remark \ref{rk100}, has $d_1=m=k+1$  and $N=d_1=m$, since the union of the $k$ lines through $q=(1:0:0)$ and the  line $x=0$  contains all the multiple points of $\A_k^x$. Take for this arrangement $p=(0:0:1)$ and note that
 $m_p=m=k+1=d_1$. The only point not connected to $p$ by lines in
 $\A_k^x$ is the point $q=(1:0:0)$. The line $L:y+2z=0$ contains $q$ and
 corresponds to Case 1 in the proof of Theorem \ref{thm1000}.
 In this case
 $$d-2m=(3k+1)-2(k+1)=k-1=m_q-1,$$
 confirming that  $\A_k^x$ is a free arrangement.
 The line $L': z=0$  contains $q$ and
 corresponds to Case 2 in the proof of Theorem \ref{thm1000}.
 In this case we have again
 $$d-2m=(3k+1)-2(k+1)=k-1=m_q-1,$$
 and this confirms again that  $\A_k^x$ is a free arrangement.
 Finally, the line $L'':y-z=0$  contains $q$ and
 corresponds to Case 3 in the proof of Theorem \ref{thm1000}.
 In this case
 $$d-2m=(3k+1)-2(k+1)=k-1=m_q-1,$$
 since all the other multiple points on $L''$ are connected to $p$ by lines in $\A_k^x$.
This confirms in a third way that  $\A_k^x$ is a free arrangement.
\end{ex}

\end{document}